\documentclass[pdflatex]{article}

\usepackage{enumerate}
\usepackage{tikz}
\usepackage{subfigure}
\usepackage[english]{babel}
\usepackage{graphicx}
\usepackage{cite}

\usepackage{amsfonts,amssymb,amsmath,latexsym,amsthm}
\usepackage{multirow}

\newtheorem{theorem}{Theorem}[section]

\newtheorem{lemma}[theorem]{Lemma}

\newtheorem{conjecture}[theorem]{Conjecture}

\theoremstyle{definition}

\begin{document}
	
\title{Multicolor Ramsey numbers on stars versus path
	\thanks{This work was supported by the National Natural Science Foundation of China (No. 12071453), the National Key R and D Program of China (2020YFA0713100),   and the Innovation Program for Quantum Science and Technology (2021ZD0302904).}
}
\author{Xuejun Zhang$^a$,\,\, Xinmin Hou$^{a,b,c,d}$
	\\
\small
 $^a$ 
 School of Data Science\\
\small University of Science and Technology of China, Hefei, Anhui 230026, China.\\
\small $^b$School of Mathematical Sciences\\
\small University of Science and Technology of China, Hefei, Anhui 230026, China.\\
\small $^{c}$CAS Key Laboratory of Wu Wen-Tsun Mathematics\\
\small University of Science and Technology of China, Hefei, Anhui 230026, China.\\
\small $^d$Hefei National Laboratory\\
\small University of Science and Technology of China, Hefei 230088, Anhui, China\\
}
\date{}
\maketitle

\begin{abstract}
For given simple graphs $H_1,H_2,\dots,H_c$, the multicolor Ramsey number $R(H_1,H_2,\dots,H_c)$ is defined as the smallest positive integer $n$ such that for an arbitrary edge-decomposition $\{G_i\}^c_{i=1}$ of the complete graph $K_n$, at least one $G_i$ has a subgraph isomorphic to $H_i$. Let $m,n_1,n_2,\dots,n_c$ be positive integers and $\Sigma=\sum_{i=1}^{c}(n_i-1)$. Some bounds and exact values of  $R(K_{1,n_1},\dots,K_{1,n_c},P_m)$ have been obtained in literature. 
Wang (Graphs Combin., 2020) conjectured that if $\Sigma\not\equiv 0\pmod{m-1}$ and $\Sigma+1\ge (m-3)^2$, then $R(K_{1,n_1},\ldots, K_{1,n_c}, P_m)=\Sigma+m-1.$ In this note, we give a new lower bound  and some exact values of $R(K_{1,n_1},\dots,K_{1,n_c},P_m)$ 
when $m\leq\Sigma$, $\Sigma\equiv k\pmod{m-1}$, and $2\leq k \leq m-2$. 
These results partially confirm Wang's conjecture. 

\end{abstract}


\section{Introduction}
In this paper,  all graphs are finite and simple.  The vertex set, edge set, minimum degree, maximum degree, complement graph and edge chromatic number of a graph $G$ are denoted by $V(G)$, $E(G)$, $\delta(G)$, $\Delta(G)$, $\overline{G}$ and $\chi'(G)$, respectively. For $v\in V(G)$, we use $d_G(v)$ to denote the degree of $v$ in $G$. We use $kG$ to denote $k$ vertex disjoint copies of $G$.   An decomposition of a graph $G$ is a set $\{G_1,G_2,\dots,G_c\}$ of edge-disjoint subgraphs of G such that $\bigcup_{i=1}^cG_i=G$. We use $K_n^r$ to denote a complete $n$-partite graph with all partition classes having the same size $r$.

Let $H_1,H_2,\dots,H_c$ be given simple graphs. The {\it multicolor Ramsey number} $R(H_1,H_2,\dots,H_c)$ is defined as the smallest positive integer $n$ such that for an arbitrary edge-decomposition $\{G_i\}^c_{i=1}$ of the complete graph $K_n$, at least one $G_i$ has a subgraph isomorphic to $H_i$.

There are few known results about $R(H_1,H_2,\ldots,H_c)$ for $c\geq 3 $ even for very special graphs. In this note, we focus on 
$R(K_{1,n_1},K_{1,n_2},\dots,K_{1,n_c},P_m)$, where $P_m$ is a path of order $m$. 
In the rest of this paper, let $\Sigma = \sum_{i=1}^{c}(n_i-1)$. Some related results are listed here.

\begin{itemize}
\item[(a)]\label{THM: grneral-upper}	(Burr, Roberts\cite{RNS}, 1973)
\begin{equation}
	R(K_{1,n_1},K_{1,n_2},\dots,K_{1,n_c})=
	\begin{cases}
		\Sigma+1, &\text{if $\Sigma$ and at least one $n_i$ are even,}\\
		\Sigma, &\text{else.}
	\end{cases}
\end{equation}

\item[(b)]\label{Upperbound}(Zhang and Zhang~\cite{strn}, 1995)
	$R(K_{1,n_1},\dots,K_{1,n_c},P_m)\leq m+\Sigma$. Moreover, if $\Sigma \equiv 0\pmod {m-1}$ and $m+\Sigma$ is odd, then $$R(K_{1,n_1},\dots,K_{1,n_c},P_m)=\Sigma+m.$$

\item[(c)]\label{origin} (Wang~\cite{smrn}, 2020)
\begin{itemize}
\item[(c1)]
If $m\geq\Sigma+1$, then $R(K_{1,n_1},\ldots,K_{1,n_c},P_m)=\max\{m,2\Sigma+1\}$.
		
\item[(c2)] If $m\leq\Sigma$, $\Sigma\equiv0 \pmod{m-1}$ and $m+\Sigma$ is even, then		
		$$R(K_{1,n_1},\dots,K_{1,n_c},P_m)=
			\begin{cases}
				\Sigma+m, &\text{if all $n_i(1\leq i \leq c)$ are odd;}\\
				\Sigma+m-1, &\text{if some $n_i$ is even.}
			\end{cases}$$
		
\item[(c3)] If $m\leq\Sigma$ and $\Sigma\not\equiv0\pmod{m-1}$, then $$R(K_{1,n_1},\dots,K_{1,n_c},P_m)\leq\Sigma+m-1.$$ Moreover, equality holds if $\Sigma\equiv1\pmod{m-1}$.

\item[(c4)] If $m\leq\Sigma$, $\Sigma\equiv k\pmod{m-1}$ and $2\leq k \leq m-2$, then $$R(K_{1,n_1},\dots,K_{1,n_c},P_m)\geq\Sigma+k-2$$.
\end{itemize}

\end{itemize}

In~\cite{smrn}, the author also proposed the following conjecture.
\begin{conjecture}[Wang~\cite{smrn}]\label{CONJ: c1}
	If $\Sigma\not\equiv 0\pmod{m-1}$ and $\Sigma+1\ge (m-3)^2$, then $$R(K_{1,n_1},\ldots, K_{1,n_c}, P_m)=\Sigma+m-1.$$
\end{conjecture}	

In this paper, we continue the study of the Ramsey number on stars versus path and give a new lower bound  and some exact values of $R(K_{1,n_1},\dots,K_{1,n_c},P_m)$
when $m\leq\Sigma$, $\Sigma\equiv k\pmod{m-1}$, and $2\leq k \leq m-2$. These results partially confirm Conjecture~\ref{CONJ: c1}. The main result is as the following.
\begin{theorem}\label{THM: main}
Let $m\leq\Sigma$, $\Sigma\equiv k\pmod{m-1}$, $2\leq k \leq m-2$, and let $s =\frac{\Sigma-k}{m-1}$.	Then the following holds. 
	\begin{itemize}
		\item[(1).]  If $s + k \geq m - 2$ and $m + \Sigma$ is even, then $$R(K_{1,n_1},\dots,K_{1,n_c},P_m) = \Sigma + m - 1.$$
		\item[(2).] If $\Sigma \equiv 0\pmod{m - 2}$, $m + \Sigma$ is odd and $n_1,\dots,n_c$ are odd, then $$R(K_{1,n_1},\dots,K_{1,n_c},P_m) = \Sigma + m - 1.$$
		\item[(3).] If $\Sigma \equiv 0\pmod{m - 2}$, $m + \Sigma$ is odd, $\Sigma < (m - 2)^2$ and some $n_i$ is even, then $$R(K_{1,n_1},\dots,K_{1,n_c},P_m) = \Sigma + m - 2.$$
		\item[(4).] $$R(K_{1,n_1},\dots,K_{1,n_c},P_m) \geq \Sigma + \min\{s + k,m - 2\}.$$
	\end{itemize}

\end{theorem}

The rest of this note is arranged as follows. We give some preliminaries in the next section. In Section 3, we prove Theorem~\ref{THM: main}.

\section{Preliminaries}
The following lemmas will be used in the proofs of this note.
We first give three fundamental results in graph theory.
\begin{lemma}[\cite{BM76}]\label{LEM: odddegree}
Every graph has an even number of odd degree vertices.
\end{lemma}

\begin{lemma}[\cite{BM76}]\label{LEM: odddegree}
For every simple graph $G$, $\chi'(G)=\Delta(G)$ or $\Delta(G)+1$.
\end{lemma}

Let $p(G)$ be the order of a longest path of a graph $G$.

\begin{lemma}[Dirac \cite{stag}]\label{path}
	Let $G$ be a connected graph of order $n\geq3$ with $\delta(G) = \delta$. Then $p(G)\geq \min\{2\delta+1,n\}$.
\end{lemma}

A graph $G$ of order $n$ is called {\it overfull} if $\lvert E(G)\rvert > \lfloor \frac{n}{2}\rfloor\Delta(G)$. Obviously, the order of an overfull graph must be odd.
\begin{lemma}[Hoffman, Rodger \cite{hg}]\label{cmg}
 Let $G$ be a complete multipartite graph. If G is not overfull, then $\chi'(G)=\Delta(G)$.	
\end{lemma}

\begin{lemma}[Omidi, Raeisi, Rahimi \cite{svsrn}]\label{chi}
	Let $n_1,\dots,n_c$ be positive integers, $\Sigma = \sum_{i=1}^{c}(n_i-1)$ and let $H$ be a graph with $\chi'(H)\leq\Sigma$. Then H can be decomposed into edge-disjoint subgraphs $H_1,H_2,\dots,H_c$ such that $\Delta(H_i)\leq n_i-1$.	
\end{lemma}

\begin{lemma}[Auerbach, Laskar \cite{la}]\label{hc}
	For  odd $n$ and  even $r$, the complete multipartite graph $K_n^r$ is the union of $\frac{n(r-1)-1}2$ edge-disjoint Hamilton cycles and a 1-factor.
	Otherwise, $K_n^r$ is the union of $\frac {n(r-1)}2$ edge-disjoint Hamilton cycles.
\end{lemma}

\section{Proof of Theorem~\ref{THM: main}}
\begin{proof}[Proof of Theorem~\ref{THM: main}]

(1) By (c3), we have $R(K_{1,n_1},\dots,K_{1,n_c},P_m) \leq \Sigma + m - 1$. Thus we only need to prove $R(K_{1,n_1},\dots,K_{1,n_c},P_m) > \Sigma + m - 2$ in this case. 
Note that $$\Sigma + m - 2=s(m-1)+k+m-2=(s - m + k+2)(m - 1) + (m - k)(m - 2).$$
Consider a decomposition {$K_{\Sigma+k-2}=G\cup\overline{G}$, where $\overline{G}=(s-m+k+2)K_{m-1}\cup(m-k)K_{m-2}$.} Obviously, $\overline{G}$ contains no $P_m$. Since $\Sigma +m$ is even, $G$ is not overfull. Note that $G$ is a complete $(s+2)$-partite graph. By Lemma~\ref{cmg}, we have $\chi^\prime(G) = \Delta(G) = \Sigma$. By Lemma~\ref{chi}, $G$ can be decomposed into edge-disjoint subgraphs $G_1,\dots,G_c$ such that $\Delta(G_i) \leq n_i - 1$. Thus $R(K_{1,n_1},\dots,K_{1,n_c},P_m) > \Sigma + m - 2$. Therefore, we have $$R(K_{1,n_1},\dots,K_{1,n_c},P_m) = \Sigma + m - 1.$$

(2) With the same reason as in (2), we need only to show $R(K_{1,n_1},\dots,K_{1,n_c},P_m) > \Sigma + m - 2$. As $n_1,\dots,n_c$ are odd, $\Sigma=\sum_{i=1}^c(n_i-1)$ is even. Since $\Sigma+m$ is odd, we have $m$ is odd. Let $\Sigma = t(m-2)$ for some positive integer $t$. Then $\Sigma+m-2=(t+1)(m-2)$.
Consider the complete (t+1)-partite graph $G$ on $\Sigma + m - 2$ vertices with all parts of the same size $m - 2$. 
 By Lemma~\ref{hc}, $G$ can be decomposed into $\frac{\Sigma}2$ edge-disjoint Hamilton cycles. For $1\leq i \leq c$, let $G_i$ be the union of $\frac{n_i-1}2$ edge-disjoint Hamilton cycles. Then $G_i$ is $(n_i-1)$-regular and thus it contains no $K_{1,n_i}$. Clearly, $\overline{G}$ contains no $P_m$. Therefore, $$R(K_{1,n_1},\dots,K_{1,n_c},P_m) = \Sigma + m - 1.$$

(3) Let $\Sigma = t(m-2)$. For the lower bound, consider the complete (t+1)-partite graph $G$ on $\Sigma +m - 3$ vertices with t parts of the same size $m-2$ and one part of size $m-3$. Since $\Sigma+m$ is odd, $\Sigma+m-3$ is even. So $G$ is not overfull. By Lemma~\ref{cmg}, $\chi^\prime(K)  = \Sigma$. By Lemma~\ref{chi}, $G$ can be decomposed into edge-disjoint subgraphs $G_1,\dots,G_c$ such that $\Delta(G_i) \leq n_i - 1$. Clearly, the complement of $G$ contains no $P_m$. Therefore, we have $$R(K_{1,n_1},\dots,K_{1,n_c},P_m) \geq \Sigma + m - 2.$$

For the upper bound, we assume to the contrary that there is a decomposition $K_{\Sigma+m-2} = \bigcup_{i=1}^{c+1} G_i$ such that $G_{i}$ contains no $K_{1,n_i}(1\leq i \leq c)$ and $G_{c+1}$ contains no $P_m$. Then $\delta(G_{c+1})\geq m-3$ (otherwise, at least one $G_i$ for $1\le i\le c$ has $\delta(G_i)\ge n_i$, a contradiction to the assumption that $G_i$ contains no $K_{1,n_i}$). 
Therefore, every component of $G_{c+1}$ has at least $m-2$ vertices.
By Lemma~\ref{path}, each component of $G_{c+1}$ have at most $m-1$ vertices. Since $\Sigma < (m-2)^2$, we have $\left| V(G_{c+1})\right|=\Sigma+m-2 <(m-2)(m-1)$. This forces that every component of $G_{c+1}$ has exact $m-2$ vertices. 
Since $\delta(G_{c+1})\geq m-3$, all components of $G_{c+1}$ are complete graphs, i.e. $G_{c+1}=\cup_{i=1}^{t+1}K_{m-2}$. Thus $d_{\cup_{i=1}^{c} G_i}(v) = \Sigma+m-3-d_{G_{c+1}}(v)=\Sigma=\sum_{i=1}^c(n_i-1)$ for all $v\in V(K_{\Sigma+m-2})$. As each $G_i$ contains no $K_{1,n_i}$, we have every $G_i$ is a $(n_i-1)$-regular graph for $1\leq i\leq c$. Note that $V(G_i)=\Sigma+m-2$ is odd. Lemma~\ref{LEM: odddegree} forces that all $n_i-1$ are even. This contradicts to the assumption that at least one $n_i$ is even.

(4) If $s + k < m - 2$, let $n=\Sigma+s+k$, then \begin{eqnarray*}
	n-1&=&\Sigma+s+k-1\\
	&=&s(m-1)+k+s+k-1\\
	&=&(s - 1)(m - 2) + m - 3 + 2(s + k).
\end{eqnarray*}
Consider a complete $(s + 2)$-partite graph $K$ on $n- 1$ vertices with $s - 1$ parts of the same size $m - 2$, one part of size $m - 3$ and two parts of size $s + k$. By Lemma~\ref{LEM: odddegree}, $\chi^\prime(K) \leq \Delta(K) + 1 = \Sigma$. By Lemma~\ref{chi}, $K$ can be decomposed into edge-disjoint subgraphs $G_1,\dots,G_c$ such that $\Delta(G_i) \leq n_i - 1$. Obviously, the complement of $K$ contains no $P_m$. So $R(K_{1,n_1},\dots,K_{1,n_c},P_m) \geq \Sigma + s + k.$

If $s + k > m - 2$, let $n=\Sigma + m - 2$, then 
\begin{eqnarray*}
	n-1&=&\Sigma + m - 3\\
	&=&s(m-1)+k+m-3\\
	&=&(s - m+k + 1)(m - 1) + (m-k + 1)(m - 2).
\end{eqnarray*}
Consider a complete $(s + 2)$-partite graph $K$ on $n-1$ vertices with $s-m+k+1$ parts of the same size $m - 1$ and $m-k+1$ parts of the same size $m - 2$. By Lemma~\ref{LEM: odddegree}, $\chi^\prime(K) \leq \Delta(K) + 1 = \Sigma$. By Lemma~\ref{chi}, $K$ can be decomposed into edge-disjoint subgraphs $G_1,\dots,G_c$ such that $\Delta(G_i) \leq n_i - 1$. Since the complement of $K$ contains no $P_m$, we have $R(K_{1,n_1},\ldots,K_{1,n_c},P_m) \geq \Sigma + m - 2.$

If $s+k=m-2$ and $m+\Sigma$ is even, then by  (1), we have $$R(K_{1,n_1},\dots,K_{1,n_c},P_m) = \Sigma + m - 1.$$

If $s+k=m-2$ and $m+\Sigma$ is odd, $\Sigma = s(m-1)+k = s(m-2)+s+k = (s+1)(m-2)< (m-2)^2$. By (2) and (3), we have
$$R(K_{1,n_1},\dots,K_{1,n_c},P_m)=
\begin{cases}
	\Sigma+m-1 &\text{if $n_i(1\leq i \leq c)$ are odd;}\\
	\Sigma+m-2 &\text{if some $n_i$ is even.}
\end{cases}$$

The proof of the theorem is completed.

\end{proof}

\section{Remarks}
In this note, we confirm some special cases ((1) and (2) in Theorem~\ref{THM: main}) of Conjecture~\ref{CONJ: c1}, and result (3) also supports the conjecture, but unfortunately, Conjecture~\ref{CONJ: c1} is still open.




\end{document}